\documentclass[12pt,reqno]{amsart}
\usepackage{epsfig}
\usepackage{graphicx}
\usepackage{mathrsfs}
\usepackage{verbatim}
\usepackage{amssymb}
\usepackage{latexsym}
\usepackage{amsmath}
\usepackage{amsfonts}
\usepackage{amsbsy}
\usepackage{amscd}
\usepackage[usenames]{color} 

\newtheoremstyle{compactcitation}{6pt}{6pt}{\itshape}{}%
{\bfseries}{.}{.5em}{\thmname{#1}\thmnumber{ #2}\thmnote{(#3)}}
\theoremstyle{compactcitation}
\newtheorem{prop}{Proposition}[section]
\newtheorem{lem}[prop]{Lemma}

\newtheorem{thm}[prop]{Theorem}
\newtheorem{conj}{Conjecture}[section]

\usepackage[T1]{fontenc}
\usepackage{lmodern}
\usepackage[noadjust]{cite}
\usepackage[hidelinks]{hyperref}
\hypersetup{pdftitle={Integer reciprocal contraction ratios of self-similar spectral measures}}
\numberwithin{equation}{section}
\newcommand{\ZZ}{\mathcal Z}
\newcommand{\supp}{\operatorname{supp}}

\newif\ifreview
\reviewtrue
\definecolor{editgreen}{rgb}{0,0.42,0.12}

\begin{document}
\baselineskip=17pt
\title[Integer reciprocal contraction ratios]
{Integer reciprocal contraction ratios of self-similar spectral measures}
\author{Zheng-Yi Lu}
\address{School of Mathematics, Hunan University, Changsha, Hunan 410082, China}
\email{zyluhnsd@163.com}
\thanks{The research is supported by NSFC of China (Nos. 12501111, 11831007, 12071125), the Hunan Provincial NSF (No. 2025JJ60065), and the
Fundamental Research Funds for the Central Universities.}
\begin{abstract}
Let $0<\rho<1$ and let $D\subset\mathbb Z$ be finite with at least two
elements. We prove that spectrality of the equal-weight self-similar
measure $\mu_{\rho,D}$ implies $\rho^{-1}\in\mathbb N$, establishing
the contraction-ratio assertion of the \L aba--Wang conjecture for integer
digit sets. We also show that an infinite orthogonal family of
exponentials forces $\rho^r\in\mathbb Q$ for some positive integer $r$.
The proof uses a finite partition of orthogonal sets into additive
cosets of $\mathbb Q(\rho)$, the rational-or-transcendental dichotomy for
mask zeros, and completeness obstructions for convolution measures.
\end{abstract}
\maketitle
\begingroup\setcounter{tocdepth}{2}\setlength{\parskip}{0pt}\tableofcontents\endgroup

\section{Introduction}\label{sect1}
\subsection{Spectral sets and singular spectral measures}
A compactly supported Borel probability measure $\mu$ on $\mathbb R^d$
is called a \emph{spectral measure} if there is a countable set
$\Lambda\subset\mathbb R^d$ such that
\[
 E(\Lambda)=\{e_\lambda(x)=e^{2\pi\mathrm{i}\langle\lambda,x\rangle}:
 \lambda\in\Lambda\}
\]
is an orthonormal basis of $L^2(\mu)$. The set $\Lambda$ is a
\emph{spectrum} of $\mu$. If only orthogonality is assumed, we call
$\Lambda$ an orthogonal set for $\mu$.

The spectral-set problem originates in Fuglede's work on commuting
self-adjoint differential operators \cite{F74}. For a measurable set
$\Omega\subset\mathbb R^d$ with $0<|\Omega|<\infty$, spectrality of
$|\Omega|^{-1}\mathbf1_\Omega\,dx$ is equivalent to the existence of
an orthogonal exponential basis for $L^2(\Omega)$. Fuglede conjectured
that this occurs precisely when $\Omega$ tiles $\mathbb R^d$ by
translations, meaning that some discrete set $T$ satisfies
\[
 \sum_{t\in T}\mathbf1_\Omega(x-t)=1\quad\text{for almost every }x.
\]
The lattice case links a tiling lattice with its dual frequency lattice
\cite{F74}. The general problem connects Fourier orthogonality with
geometric decompositions and arithmetic properties of translation sets.

Tao \cite{T04} disproved the unrestricted conjecture in dimensions at
least five. 
Matolcsi \cite{M05} reduced
the dimension to four, and Kolountzakis and Matolcsi
\cite{KM06} subsequently constructed a counterexample
in dimension three.
Substantial positive results have also been obtained. Lev and Matolcsi
\cite{LM22} proved the conjecture for convex bodies in every dimension,
while \L aba and Londner \cite{LL23} established the Coven--Meyerowitz
tiling conditions for a class of integer tilings involving three
distinct odd prime factors. 
More recently, Zhang \cite{Z26} announced planar counterexamples to both implications.
Despite these counterexamples, the connections between spectrality and
tiling continue to attract substantial research, particularly for
geometrically restricted sets and arithmetic classes of tilings
\cite{L01, LM22, LL23}.

The measure formulation also includes singular probability measures,
for which Lebesgue measure of the support may vanish. Jorgensen and
Pedersen \cite{JP98} showed that the middle-fourth Cantor measure is
spectral, whereas the middle-third Cantor measure is not. Thus fractal
support is compatible with a Fourier basis, but self-similarity alone
does not imply spectrality. The resulting theory studies Fourier
expansions in $L^2(\mu)$ on fractal probability spaces and reveals
arithmetic phenomena absent from the classical interval example.

	Since the work of Jorgensen and Pedersen, singular spectral measures
have attracted sustained and growing interest among researchers in
harmonic analysis and fractal geometry. The literature now extends from
self-similar and self-affine measures to nonstationary and random
infinite convolutions \cite{AL23, LW02,LDL22, DHL19,AFL19,LMW22}.
Construction methods often start with finite Fourier matrices and
iterate them across scales. In one dimension, \L aba and Wang
\cite{LW02} established spectrality under a compatible-pair condition,
and Dutkay, Haussermann, and Lai \cite{DHL19} proved that integer
Hadamard triples generate self-affine spectral measures in arbitrary
dimension.

Beyond the existence of exponential bases, the theory investigates
convergence of Fourier expansions, dilation properties of spectra, and
their quantitative distribution. Strichartz \cite{S00,S06} studied
convergence of mock Fourier series, while Dutkay, Han, and Sun
\cite{DHS14} obtained divergence results showing that convergence can
depend on the chosen spectrum. Spectral eigenvalue problems ask which
real numbers $t$ allow both $\Lambda$ and $t\Lambda$ to be spectra of
the same measure; related structure and existence results appear in
\cite{HTW19,KLW26,L26}. Beurling dimension describes the large-scale
counting growth of frequency sets: dimension bounds were developed in
\cite{DHSW11,HKTW18, ILLW22, S21}, and intermediate-value results for spectra were
obtained in \cite{LW24,L26}. Together these developments demonstrate the
variety of Fourier bases available for singular measures and the
interaction between their arithmetic and geometric properties.
The present paper concerns a basic existence question: which
contraction ratios can support an exponential orthonormal basis?

\subsection{Self-similar measures and the \L aba--Wang conjecture}
Fix $0<\rho<1$ and a finite set $D\subset\mathbb Z$ of cardinality
$N\geq2$. Hutchinson's theorem \cite{H81} gives a unique compactly
supported probability measure satisfying
\begin{equation}\label{(1.1)}
 \mu_{\rho,D}(E)=\frac1N\sum_{d\in D}
 \mu_{\rho,D}(\rho^{-1}E-d)
 \qquad(E\subset\mathbb R\text{ Borel}).
\end{equation}
The associated contractions are $x\mapsto\rho(x+d)$.
We recall the original conjectural formulation of \L aba and Wang,
restricted here to positive contractions.
\begin{conj}[{\cite[Conjecture~3.1]{LW02}}]\label{conj1.1}
Let $A=\{a_1,\ldots,a_N\}\subset\mathbb R$ consist of $N\geq2$
distinct digits, with $0\in A$, and let $p_j>0$ with $\sum_jp_j=1$.
Suppose that the self-similar measure associated with
$x\mapsto\rho(x+a_j)$, $0<\rho<1$, and weights $p_j$ is spectral.
Then:
\begin{enumerate}
\item $\rho^{-1}=q$ is an integer with $q\geq2$;
\item $p_1=\cdots=p_N=1/N$;
\item There exists $\alpha>0$ such that $A=\alpha D_0$, where $D_0\subset \mathbb{Z}$ and $D_0$ tiles $\mathbb{Z}_q$, i.e., there exists $C$ such that $D_0\oplus C=\mathbb{Z}_q\pmod{q}$. 
\end{enumerate}
\end{conj}
The three clauses above
record the original conjecture; our result concerns its first clause
under the hypothesis $D\subset\mathbb Z$, and does not assert the
third clause. The equal-weight assumption is supported by the uniformity
result of Deng and Chen \cite{DC21}.

For Bernoulli convolutions, Dai \cite{D12} proved that spectrality is
equivalent to the contraction ratio being the reciprocal of an even
integer. Dai, He, and Lau \cite{DHL14} characterized spectral
$N$-Bernoulli measures. These results establish the integer-reciprocal
restriction for important classes of digit sets. 
For more general digit sets, An and Wang \cite[Theorem~1.3]{AW21}
proved the integer reciprocal-contraction conclusion when
\begin{equation}\label{(1.2)}
 \ZZ(\widehat{\delta_D})\subset\alpha\mathbb Z
 \quad\text{for some }\alpha>0,
 \qquad \delta_D=\frac1N\sum_{d\in D}\delta_d,
\end{equation}
where $\ZZ(f)=\{\xi\in\mathbb R:f(\xi)=0\}$.
Although integer digits make $\widehat{\delta_D}$ periodic with only
finitely many zeros modulo one, those zeros need not all be rational.
Consequently, periodicity alone does not imply the lattice assumption
in \eqref{(1.2)}. Our purpose is to remove that assumption for integer
digit sets.

\subsection{Main results and proof strategy}
Our principal result is the following.
\begin{thm}\label{thm1.1}
Let $0<\rho<1$ and let $D\subset\mathbb Z$ be finite with $\#D\geq2$.
If $\mu_{\rho,D}$ is a spectral measure, then
$
 \rho^{-1}\in\mathbb N.
$
\end{thm}
Thus the contraction-ratio assertion of Conjecture~\ref{conj1.1} holds
for arbitrary finite integer digit sets with equal weights.

A related, but weaker, question asks when $L^2(\mu)$ contains an
infinite orthogonal family of exponentials, without requiring
completeness. Hu and Lau \cite{HL08} proved that a Bernoulli convolution
has such a family precisely when $\rho^r=p/q$ for some positive integer
$r$, where $p$ is odd and $q$ is even in lowest terms. Deng
\cite[Theorem~1.2]{D14} obtained the corresponding characterization for
consecutive digit sets of prime cardinality. Infinite orthogonality and
bounds on the size of orthogonal families have also been studied for
higher-dimensional self-affine measures. 
These works illustrate the continuing interest in the arithmetic
restrictions imposed by infinite orthogonality, independently of
spectral completeness \cite{HL08,D14,LDZ19}.
Our second theorem provides a rational-power restriction for every
finite integer digit set and is the main arithmetic step in the proof
of Theorem~\ref{thm1.1}.

\begin{thm}\label{thm1.2}
Let $0<\rho<1$ and let $D\subset\mathbb Z$ be finite with $\#D\geq2$.
If $L^2(\mu_{\rho,D})$ contains infinitely many mutually orthogonal
exponential functions, then $\rho^r\in\mathbb Q$ for some integer
$r\geq1$.
\end{thm}
Theorem~\ref{thm1.2} extends the rational-power restriction familiar
from the Bernoulli case \cite{HL08} to arbitrary finite integer digit
sets. It does not imply spectrality: Theorem~\ref{thm1.1} additionally
uses completeness to exclude irrational roots of rational numbers and
rational contractions with numerator greater than one.

Set $b=\rho^{-1}$ and $K=\mathbb Q(b)$. The main structural step,
Lemma~\ref{thm3.1}, partitions every orthogonal set into finitely many
additive cosets of $K$, with differences between distinct classes
confined to finitely many Fourier-zero levels. The elementary equation
$b^M-Ab^R=1$ controls the possible differences of levels outside $K$.
Within a coset, a polynomial argument excludes transcendental $b$.
The Gelfond--Schneider theorem then implies that every relevant
algebraic mask zero is rational. Ramsey's theorem and Deng's
prime-digit condition \cite{D14} complete the rational-power reduction.
For a spectrum, an endpoint Bessel estimate excludes the irrational-root
case. The rational case is treated using a common denominator bound,
Gao--Ma Fourier decay \cite{GM17}, and the finite-section criterion
in \cite[Theorem~2.5(ii)]{AW21}.

Section~\ref{sect2} collects the common Fourier and algebraic facts.
Section~\ref{sect3} proves the finite partition lemma and
Theorem~\ref{thm1.2}. Section~\ref{sect4} proves
Theorem~\ref{thm1.1}, introducing the remaining auxiliary results at
their points of use.

\section{Preliminary lemmas}\label{sect2}
\subsection{Fourier zeros and orthogonality}
It is useful to view the measure probabilistically. Let $X_1,X_2,\ldots$
be independent random variables, each uniformly distributed on $D$.
Then $\mu$ is the distribution of the absolutely convergent random series
$\sum_{j\geq1}\rho^jX_j$. The measure $\mu_n$ is the distribution of
its first $n$ terms. The remaining terms have the same distribution as
$\rho^n$ times an independent copy of the full series. This explains the
convolution identity below. Here convolution is defined by
\[
 \int f\,d(\nu*\eta)=\iint f(x+y)\,d\nu(x)\,d\eta(y).
\]
A finite convolution $\mu_n$ has at most $N^n$ support points, whether or
not different digit sequences produce the same point.
For a probability measure $\nu$, we use
\[
 \widehat\nu(\xi)=\int e^{-2\pi\mathrm{i}\xi x}\,d\nu(x),
 \qquad \delta_D=\frac1N\sum_{d\in D}\delta_d.
\]
Then $\mu:=\mu_{\rho,D}=\mathop{*}_{n=1}^\infty\delta_{\rho^nD}
=\mu_n*\mu_{>n}$, where 
\begin{equation}\label{(2.1)}
\mu_n=\mathop{*}_{j=1}^n\delta_{\rho^jD},
 \qquad \mu_{>n}=\mathop{*}_{j=1}^\infty\delta_{\rho^{n+j}D}.
\end{equation}
Let the digit mask be 
\begin{align*}
	m_D(\xi)=\widehat{\delta_D}(\xi)=\frac1N\sum_{d\in D}e^{-2\pi\mathrm{i}d\xi}.
\end{align*}
Since $D\subset\mathbb Z$, there is a finite set $Z\subset(0,1)$ such that
\begin{equation}\label{(2.2)}
 \ZZ(m_D)=Z+\mathbb Z,\qquad
 \ZZ(\widehat\mu)=\bigcup_{n\geq1}\mathcal Z_n,
 \qquad\mathcal Z_n=b^n(Z+\mathbb Z),\quad b=\rho^{-1}.
\end{equation}
The zero-set notation in \eqref{(2.2)} follows from the periodicity of
$m_D$ and the product $\widehat\mu(\xi)=\prod_{n\geq1}m_D(\rho^n\xi)$.
A set $\Lambda$ is orthogonal for $\mu$ if and only if
\begin{equation}\label{(2.3)}
	(\Lambda-\Lambda)\setminus\{0\}\subset\ZZ(\widehat\mu).
\end{equation}
We use the following classical criterion without proof.
\begin{lem}[{\cite{JP98}}]\label{thm2.1}
Let $\mu$ be a compactly supported Borel probability measure on
$\mathbb R$, and let $\Lambda\subset\mathbb R$ be countable. Then
$\Lambda$ is orthogonal for $\mu$ if and only if
\begin{equation}\label{(2.4)}
 Q_{\mu,\Lambda}(\xi):=\sum_{\lambda\in\Lambda}
 |\widehat\mu(\xi+\lambda)|^2\leq1,\;\;\forall\;\xi\in\Bbb R. 
\end{equation}
Furthermore, $\Lambda$ is a spectrum if and only if equality holds for every real
$\xi$.
\end{lem}

Translating an orthogonal set preserves orthogonality and completeness.
We shall therefore assume $0\in\Lambda$ whenever convenient. Translating
$D$ also only translates the resulting measure; in endpoint arguments
we normalize $\min D=0$.

\subsection{Algebraic preliminaries}
The first lemma concerns ordinary integer powers, including negative
powers. Its proof uses only the fact that the integers have no finite
accumulation point. Both exponents are unknown, whereas $b$ and $A$
are fixed.
\begin{lem}\label{thm2.2}
Let $b>1$ and $A\in\mathbb R\setminus\{0\}$. The equation
\begin{equation}\label{(2.5)}
 b^M-Ab^R=1,\qquad (M,R)\in\mathbb Z^2,
\end{equation}
has only finitely many solutions.
\end{lem}
\begin{proof}
A sequence of integer powers of $b$ converging to a nonzero finite limit
is eventually constant: taking logarithms gives a convergent sequence of
integers. For a fixed $M$, there is at most one $R$ in \eqref{(2.5)}.
If infinitely many solutions existed, we could select distinct $M_k$
tending to $+\infty$ or $-\infty$. In the latter case,
$b^{R_k}=(b^{M_k}-1)/A\to-1/A$. This is impossible if the limit is
negative, and otherwise forces $R_k$, and hence $M_k$, to be eventually
constant. In the former case,
\[
 b^{M_k-R_k}=\frac{A}{1-b^{-M_k}}\longrightarrow A.
\]
Again either the limit is impossible or $M_k-R_k=h$ eventually and
$b^h=A$. Then $b^{M_k}-Ab^{R_k}=0$, a contradiction.
\end{proof}

The partition lemma in Section~\ref{sect3} uses a general real subfield.
In its applications
we take $K=\mathbb Q(b)$. The quotient $\mathbb R/K$ is used only as a
vector space over $K$; no quotient-field structure is intended.
For clarity, $\mathbb Q(b)$ is the smallest field containing the rational
numbers and $b$. Its elements are quotients $P(b)/Q(b)$, where
$P,Q\in\mathbb Q[x]$ and $Q(b)\ne0$. An additive coset is a set
$\gamma+K=\{\gamma+k:k\in K\}$. Thus
\[
 x\equiv y\pmod K\quad\Longleftrightarrow\quad x-y\in K.
\]
Congruences may be added and multiplied by elements of $K$. For example,
$\sqrt3+\mathbb Q(\sqrt2)$ consists of the numbers
$\sqrt3+a+c\sqrt2$ with $a,c\in\mathbb Q$. Two numbers are in the same
class precisely when their difference is in $K$.

There are two useful cases to keep in mind. If $b$ is algebraic of
degree $d$, then every element of $\mathbb Q(b)$ is a unique rational
linear combination of $1,b,\ldots,b^{d-1}$. This follows by dividing
polynomials by the minimal polynomial; inverses are obtained from the
Euclidean algorithm because that polynomial is irreducible. Thus
$\mathbb Q(\sqrt2)=\{a+c\sqrt2:a,c\in\mathbb Q\}$.
If $b$ is transcendental, no nonzero rational polynomial vanishes at
$b$, so rational functions in $b$ obey exactly the same identities as
rational functions in an indeterminate $x$.

Only elementary field operations are needed for congruences here.
If $x-y\in K$ and $a\in K$, then $a(x-y)\in K$. If also $a\ne0$,
then $ax-ay\in K$ implies $x-y\in K$ by division by $a$.
One never divides by a coset or treats $\mathbb R/K$ as a field.

An algebraic number is a root of a nonzero polynomial with rational
coefficients. A transcendental number is one that is not algebraic.
The following classical theorem is the only transcendence theorem used
in this paper; we quote it without proof.
\begin{thm}[{\cite{B75}}]\label{thm2.3}
Let $\alpha,\beta\in\mathbb C$ be algebraic numbers with
$\alpha\ne0,1$ and $\beta\notin\mathbb Q$. Then every value of
$\alpha^\beta$ is transcendental.
\end{thm}
For complex numbers, a value of $\alpha^\beta$ means
$\exp(\beta L)$ for a choice of logarithm $L$ with $\exp L=\alpha$.
We need only $\alpha=-1$ and $L=\pi\mathrm{i}$.

\begin{lem}\label{thm2.4}
Let $D\subset\mathbb Z$ be finite and nonempty. Every real zero of
$m_D$ is either rational or transcendental. 
\end{lem}
\begin{proof}
If $m_D(z)=0$, then $\alpha=e^{-2\pi\mathrm{i}z}$ is algebraic, since
it is a root of the polynomial
$\sum_{d\in D}X^{d-\min D}$. If $z$ were algebraic and irrational,
Theorem~\ref{thm2.3} would make every value of
$(-1)^{-2z}$ transcendental. Taking the logarithm $\log(-1)=\pi\mathrm{i}$
gives the value $e^{-2\pi\mathrm{i}z}=\alpha$, a contradiction.
\end{proof}
The transcendence theorem in this proof is an external input. In
particular, the fact that $\alpha$ is algebraic does not by itself mean
that it is a root of unity; that shortcut would be invalid.

\begin{lem}\label{thm2.5}
Let $b>1$, and let $r$ be the least positive integer with $b^r\in\mathbb Q$.
Then the minimal polynomial of $b$ over $\mathbb Q$ is $X^r-b^r$.
In particular, $1,b,\ldots,b^{r-1}$ are linearly independent over
$\mathbb Q$.
\end{lem}
\begin{proof}
Let $f\in\mathbb Q[X]$ be the monic minimal polynomial, of degree $d$.
Since $f$ divides $X^r-b^r$, we have $d\leq r$, and every root of $f$
has complex absolute value $b$. Its nonreal roots occur in conjugate
pairs, each with product $b^2$; every real root equals $b$ or $-b$.
Their total product is therefore $\pm b^d$. By Vieta's formula this
product is rational and nonzero, so $b^d\in\mathbb Q$. Minimality of
$r$ gives $d\geq r$. Thus $d=r$ and $f=X^r-b^r$. A nonzero rational
linear relation among $1,b,\ldots,b^{r-1}$ would contradict the
minimal degree of $f$.
\end{proof}

\section{Proof of Theorem \ref{thm1.2}}\label{sect3}
This section proves the arithmetic reduction. We first establish the
finite partition, then exclude the transcendental case and reduce
the rational zero sets to a prime-digit problem.
\subsection{Finite partitions of orthogonal sets}
The hypothesis in the next lemma is a condition on the \emph{whole}
difference set. It does not say that the differences lie in one branch
$b^n(z+\mathbb Z)$ or even in the union over $n$ for one fixed $z$.
The conclusion is a finite partition, not such a single-branch assertion.
\begin{lem}\label{thm3.1}
Let $b>1$, let $K\subset\mathbb R$ be a subfield containing $b$, and let
$Z\subset(0,1)$ be finite. Put
\[
 \mathcal S=\bigcup_{n\geq1}b^n(Z+\mathbb Z).
\]
Suppose $0\in\Lambda$ and
$(\Lambda-\Lambda)\setminus\{0\}\subset\mathcal S$.
Then $\Lambda$ meets only finitely many additive cosets of $K$.
Writing the nonempty classes as $\Lambda_1,\ldots,\Lambda_s$, there is a
finite set $F\subset\mathbb N$ such that
\begin{equation}\label{(3.1)}
 \Lambda_i-\Lambda_j\subset\bigcup_{n\in F}b^n(Z+\mathbb Z)
 \qquad(i\ne j).
\end{equation}
In addition, all representations of points in $\Lambda\setminus K$
use only finitely many pairs $(n,z)$ with $z\in Z\setminus K$.
The number of such pairs has a bound depending only on $b,K,Z$.
\end{lem}
\begin{proof}
For $z,w\in Z\setminus K$, there is at most one $c\in K$ satisfying
$w-cz\in K$. Indeed, two distinct choices imply $z\in K$.
To spell this out, if $w-c_1z=k_1$ and $w-c_2z=k_2$ with
$k_1,k_2\in K$ and $c_1\ne c_2$, subtraction gives
$(c_2-c_1)z=k_1-k_2$. Dividing by $c_2-c_1\in K\setminus\{0\}$
would give $z\in K$, a contradiction. Whenever
such a $c$ exists, it is nonzero; denote it by $c(z,w)$.

Fix $z\in Z\setminus K$ and let
\[
 J_z=\{n\geq1:\Lambda\cap b^n(z+\mathbb Z)\ne\varnothing\}.
\]
For distinct $j,k\in J_z$, choose
$\lambda_j=b^j(z+u)$ and $\lambda_k=b^k(z+v)$ in $\Lambda$.
They are distinct, because otherwise $(b^j-b^k)z\in K$.
Their difference has a representation $b^\ell(w+t)$ with $w\in Z$.
Modulo $K$ this gives
\[
 (b^j-b^k)z\equiv b^\ell w\pmod K.
\]
The left side does not belong to $K$: otherwise division by the nonzero
number $b^j-b^k\in K$ would put $z$ in $K$. Thus $w\notin K$.
Dividing the congruence by $b^\ell$ shows that
\[
 w-(b^{j-\ell}-b^{k-\ell})z\in K.
\]
By uniqueness, $c(z,w)=b^{j-\ell}-b^{k-\ell}$. Multiplication by
$b^{\ell-k}$ gives
\begin{equation}\label{(3.2)}
 b^{j-k}-c(z,w)b^{\ell-k}=1.
\end{equation}
For a fixed $z$, there are only finitely many possible $w$, and hence
only finitely many possible coefficients $c(z,w)$. For each coefficient,
Lemma~\ref{thm2.2}, applied with $M=j-k$ and $R=\ell-k$, gives only
finitely many exponent pairs. In particular, there is a finite set
$E_z\subset\mathbb Z$ containing every possible nonzero difference
$j-k$. If $J_z$ is nonempty, fix $k_0\in J_z$. Every other $j\in J_z$
satisfies $j-k_0\in E_z$, and the map $j\mapsto j-k_0$ is injective.
Thus $\#J_z\leq1+\#E_z$. This bound does not depend on $\Lambda$. Summing over the finitely many $z$ proves the
last assertion. Each set $b^n(z+\mathbb Z)$ lies in one additive coset
of $K$. Together with the class $K$ itself, these give the required
finite partition.
In detail, because $0\in\Lambda$, each nonzero $\lambda\in\Lambda$
belongs to $\mathcal S$. If a representation uses $z\in K$, then
$\lambda\in K$. Otherwise it uses one of the finitely many pairs
$(n,z)$ just obtained, and
$b^n(z+\mathbb Z)\subset b^nz+K$. These finitely many cosets, together
with $K$, cover $\Lambda$. This is exactly where the assumption
$0\in\Lambda$ is used.

Choose representatives $\gamma_i\in\Lambda_i$. For $i\ne j$, set
$\delta_{ij}=\gamma_i-\gamma_j\notin K$. If
$\lambda\in\Lambda_i$ and $\lambda'\in\Lambda_j$, then
$\lambda-\lambda'\equiv\delta_{ij}\pmod K$. Every zero-set
representation $\lambda-\lambda'=b^n(w+m)$ therefore satisfies
\[
 \delta_{ij}\equiv b^nw\pmod K,\qquad w\in Z\setminus K.
\]
For fixed $i,j,w$, at most one exponent $n$ is possible. Otherwise,
$(b^n-b^{n'})w\in K$ for distinct $n,n'$. Since $b>1$, the coefficient
$b^n-b^{n'}$ is a nonzero element of $K$, and division would imply
$w\in K$, a contradiction. There are only finitely many pairs $i,j$
and finitely many $w\in Z$. Taking the union of their possible
exponents gives a finite set $F$ satisfying \eqref{(3.1)}.
In particular, the representative of the class containing zero may be
chosen to be zero.
\end{proof}

The finite partition has two distinct uses. An infinite orthogonal set
has an infinite class, by the pigeonhole principle. For a spectrum,
however, discarding all the other classes need not preserve completeness.
We therefore retain every class in Section~\ref{sect4}.

\subsection{Arithmetic reduction}
\begin{lem}\label{thm3.2}
If $b=\rho^{-1}$ is transcendental, every orthogonal set for
$\mu_{\rho,D}$ contained in an additive coset of $\mathbb Q(b)$ is
finite.
\end{lem}
\begin{proof}
The empty set is finite, so assume the orthogonal set is nonempty.
Translate it so that $0\in\Gamma\subset K:=\mathbb Q(b)$.
The assertion is immediate if $\Gamma=\{0\}$. For every nonzero
$\delta\in\Gamma-\Gamma$, orthogonality gives
$\delta=b^n(z+m)$ with $n\geq1$, $z\in Z$, and $m\in\mathbb Z$.
Since $b,\delta\in K$, we have $z=b^{-n}\delta-m\in K$.
Thus, writing $Z_K=Z\cap K$, we obtain
\[
 (\Gamma-\Gamma)\setminus\{0\}
 \subset\bigcup_{n\geq1}b^n(Z_K+\mathbb Z).
\]
In particular, $Z_K$ is nonempty when $\Gamma$ has at least two elements.

Choose a common polynomial denominator for this finite set:
$z=P_z(b)/Q(b)$, where $P_z,Q\in\mathbb Q[x]$ and $Q\ne0$.
The polynomials $P_z$ are nonzero because $0<z<1$.
Set
\[
 L=\max\bigl(\{\deg Q\}\cup\{\deg P_z:z\in Z_K\}\bigr).
\]
For each nonzero $\delta\in\Gamma-\Gamma$, a representation
$\delta=b^n(z+m)$ yields the nonzero polynomial
\[
 G_\delta(x)=x^n(P_z(x)+mQ(x)),\qquad
 G_\delta(b)=Q(b)\delta.
\]
Evaluation at the transcendental number $b$ is injective on
$\mathbb Q[x]$: two rational polynomials with the same value at $b$
must be identical. Hence $G_\delta$ is independent of its chosen
representation. Put $G_0=0$. In particular,
\[
 G_{\lambda-\lambda'}=G_\lambda-G_{\lambda'}
 \qquad(\lambda,\lambda'\in\Gamma).
\]
For every nonzero $G_\delta$, the difference between its largest and
smallest exponents with nonzero coefficients is at most $L$.
This follows directly from the factorization
$x^n(P_z+mQ)$, since $\deg(P_z+mQ)\leq L$.

Fix $\lambda_*\in\Gamma\setminus\{0\}$ and write
$d=\deg G_{\lambda_*}$. If $\deg G_\lambda>d+2L$ for some
$\lambda\in\Gamma$, then the smallest exponent in $G_\lambda$ is
greater than $d+L$. Thus its nonzero terms have exponents disjoint
from those of $G_{\lambda_*}$, all of which are at most $d$.
The polynomial $G_\lambda-G_{\lambda_*}$ retains both its highest
term, of degree greater than $d+2L$, and a nonzero term of degree
at most $d$. This contradicts the preceding bound $L$ for
$G_{\lambda-\lambda_*}$. We have therefore proved
\begin{equation}\label{(3.3)}
 \deg G_\lambda\leq H:=d+2L
 \qquad(\lambda\in\Gamma\setminus\{0\}).
\end{equation}
It follows that $\deg G_{\lambda-\lambda'}\leq H$ for distinct
$\lambda,\lambda'\in\Gamma$. In any representation
$G_{\lambda-\lambda'}=x^n(P_z+mQ)$, the second factor is a nonzero
polynomial, so $n\leq\deg G_{\lambda-\lambda'}\leq H$.
Consequently,
\[
 (\Gamma-\Gamma)\setminus\{0\}
 \subset\bigcup_{n=1}^{H}b^n(Z+\mathbb Z)
 =\ZZ(\widehat\mu_H).
\]
Thus $\Gamma$ is orthogonal for the finite measure $\mu_H$, which
has at most $N^H$ support points. The functions in $E(\Gamma)$
are linearly independent in $L^2(\mu_H)$, and hence
$\#\Gamma\leq\dim L^2(\mu_H)\leq N^H<\infty$.
\end{proof}

\begin{lem}[{\cite{R30}}]\label{thm3.3}
If the unordered pairs of a countably infinite set are colored with
finitely many colors, there is an infinite subset all of whose unordered
pairs have the same color.
\end{lem}
When this lemma is applied to differences, the sets defining the colors
below will always be symmetric under $x\mapsto-x$. If they overlap, a
fixed ordering of the colors resolves ties.

\begin{lem}[{\cite{D14}}]\label{thm3.4}
Let $b>1$ and let $\ell$ be a prime. If an infinite set
$\Gamma\subset\mathbb R$ satisfies
\[
 (\Gamma-\Gamma)\setminus\{0\}
 \subset\bigcup_{n\geq1}b^n(\mathbb Z\setminus\ell\mathbb Z),
\]
then $b^r\in\mathbb Q$ for some integer $r\geq1$.
\end{lem}
This is the necessary part of Deng's
\cite[Theorem~1.2]{D14}, written in zero-set notation. To explain the change of notation, let
$D_\ell=\{0,1,\ldots,\ell-1\}$ and $\rho=b^{-1}$. Then
\[
 \ZZ(m_{D_\ell})=\frac1\ell(\mathbb Z\setminus\ell\mathbb Z),
 \qquad
 \ZZ(\widehat\mu_{\rho,D_\ell})
 =\frac1\ell\bigcup_{n\geq1}b^n(\mathbb Z\setminus\ell\mathbb Z).
\]
Thus the displayed difference condition says exactly that
$\ell^{-1}\Gamma$ is an infinite orthogonal set for
$\mu_{\rho,D_\ell}$. Deng's necessary condition gives
$\rho^r\in\mathbb Q$, equivalently $b^r\in\mathbb Q$.
We quote that condition without reproducing its proof.

\medskip\noindent\textbf{Proof of Theorem~\ref{thm1.2}.}
If $\rho\in\mathbb Q$, the conclusion holds with $r=1$. We may
therefore assume $\rho$ is irrational and prove directly that a
positive power of $\rho$ is rational.

Let $\Lambda$ be an infinite orthogonal set for $\mu_{\rho,D}$, and
translate it so that $0\in\Lambda$. Put $b=\rho^{-1}$ and
$K=\mathbb Q(b)$. By Lemma~\ref{thm3.1}, $\Lambda$ is a finite union
of its classes modulo $K$. At least one class is infinite. Subtracting
one of its elements gives an infinite orthogonal set $\Gamma$ containing zero and satisfying 
\begin{equation}\label{(3.4)}
(\Gamma-\Gamma)\setminus\{0\}\subset K.
\end{equation}
If $b$ were transcendental, Lemma~\ref{thm3.2} would contradict
\eqref{(3.4)}. Hence $b$ is algebraic.

For distinct $\lambda,\lambda'\in\Gamma$, orthogonality and
\eqref{(3.4)} give
 \begin{align*}
	\lambda-\lambda'=b^n(z+m)\in K,
\end{align*}
which yields that $z=b^{-n}(\lambda-\lambda')-m\in K$. 
Since $b$ is algebraic, every element of $K$ is algebraic.
Lemma~\ref{thm2.4} therefore implies $z\in\mathbb Q$.
Consequently,
\begin{align*}
	(\Gamma-\Gamma)\setminus\{0\}\subset K\cap\cup_{j=1}^\infty b^j(\Bbb Q\cap Z+\Bbb Z).
\end{align*}
Let $L\geq2$ be a common denominator of the nonempty finite set
$Z\cap\mathbb Q$. Since its elements lie strictly between zero and one,
\[
 L\bigl((Z\cap\mathbb Q)+\mathbb Z\bigr)\subset\mathbb Z\setminus L\mathbb Z=\bigcup_{h=1}^t\ \bigcup_{j=0}^{e_h-1}
 \ell_h^j(\mathbb Z\setminus\ell_h\mathbb Z),
\]
where $L=\prod_{h=1}^t\ell_h^{e_h}$ and $\ell_1,\cdots,\ell_t$ are distinct primes. 
Each component on the right is symmetric. Color every unordered pair of
$L\Gamma$ by one of the finitely many pairs $(h,j)$ for which its
difference belongs to
\[
 \ell_h^j\bigcup_{n\geq1}b^n(\mathbb Z\setminus\ell_h\mathbb Z).
\]
This does not depend on the orientation of the difference.
Ramsey's theorem then
provides an infinite set $A\subset L\Gamma$ such that
\[
 (A-A)\setminus\{0\}\subset
 \ell_h^j\bigcup_{n\geq1}b^n(\mathbb Z\setminus\ell_h\mathbb Z)
\]
for one fixed pair $(h,j)$. 
Dividing that subset by $\ell_h^j$ and applying
Lemma~\ref{thm3.4} gives $b^r\in\mathbb Q$ for some $r\geq1$.
Equivalently, $\rho^r\in\mathbb Q$, proving
Theorem~\ref{thm1.2}.\hfill$\square$

\section{Proof of Theorem \ref{thm1.1}}\label{sect4}
Assume that $\mu=\mu_{\rho,D}$ is spectral and let $0\in\Lambda$ be a
spectrum. We may replace $D$ by $D-\min D$, so that
$0\in D\subset\mathbb Z_{\geq0}$. Set $b=\rho^{-1}$.
The support of $\mu$ is infinite: if $d\in D\setminus\{0\}$,
it contains the distinct points $\rho^n d$, $n\geq1$.
Thus $L^2(\mu)$ is infinite-dimensional and $\Lambda$ is infinite.
Theorem~\ref{thm1.2} now gives
\begin{equation}\label{(4.1)}
 \rho^r=p/q,\qquad 1\leq p<q,\quad\gcd(p,q)=1,
\end{equation}
where $r$ is chosen minimal. We first prove $r=1$, then $p=1$.
The two exclusions use completeness in different ways. If $r>1$, the
frequencies split into finitely many sets orthogonal for proper factors;
the endpoint estimate contradicts Parseval. If $r=1$ but $p>1$, a
summable Fourier estimate permits application of the quoted finite-section
criterion. 

\subsection{Excluding irrational roots of rational numbers}
Suppose $r>1$. By Lemma~\ref{thm2.5},
$1,b,\ldots,b^{r-1}$ are linearly independent over $\mathbb Q$.
Separate the convolution factors according to their residue classes:
\begin{equation}\label{(4.2)}
 \mu=\nu_0*\cdots*\nu_{r-1},\qquad
 \nu_s=\mathop{*}_{\substack{n\geq1\\ n\equiv s\ ({\rm mod}\ r)}}
          \delta_{\rho^nD}.
\end{equation}
Applying the infinite-product identity to this subsequence
of factors, as in \eqref{(2.2)}, gives
\[
 \ZZ(\widehat\nu_u)=
 \bigcup_{\substack{n\geq1\\n\equiv u\ ({\rm mod}\ r)}}b^n(Z+\mathbb Z).
\]
Put $K=\mathbb Q(b)$. Lemma~\ref{thm3.1} gives a finite disjoint
partition $\Lambda=\bigcup_{i=1}^s\Lambda_i$ into classes modulo $K$.
For each class choose $\gamma_i\in\Lambda_i$ and set
$\Gamma_i=\Lambda_i-\gamma_i\subset K$. 
For every nonzero difference $\delta\in\Gamma_i-\Gamma_i$, its
representation $\delta=b^n(z+m)$ has $z\in K$. By
Lemma~\ref{thm2.4}, $z$ is rational. Since $b^r=q/p\in\mathbb Q$,
\[
 (\Gamma_i-\Gamma_i)\setminus\{0\}
 \subset\bigcup_{u=0}^{r-1}b^u(\mathbb Q\setminus\{0\}).
\]
We claim that each $\Gamma_i$ lies in a single rational line
$b^{u_i}\mathbb Q$. If $\Gamma_i=\{0\}$, choose $u_i$ arbitrarily.
Otherwise choose $ab^{u_i}\in\Gamma_i\setminus\{0\}$ with
$a\in\mathbb Q\setminus\{0\}$. An element $cb^u\in\Gamma_i$
with $u\ne u_i$ and $c\ne0$ would give
\[
 cb^u-ab^{u_i}=db^w,
 \qquad d\in\mathbb Q\setminus\{0\},\quad 0\leq w<r.
\]
This contradicts the linear independence of $1,b,\ldots,b^{r-1}$:
if $w=u$, the coefficient of $b^{u_i}$ remains nonzero; if $w=u_i$,
the coefficient of $b^u$ remains nonzero; otherwise all three nonzero
coefficients belong to distinct basis elements. This proves the claim.

Every nonzero difference in $\Gamma_i$ lies in $b^{u_i}\mathbb Q$.
Writing $n=kr+v$ in a zero-set representation shows that the same
difference lies in $b^v\mathbb Q$, so linear independence forces
$v=u_i$. Hence
\[
 (\Lambda_i-\Lambda_i)\setminus\{0\}
 \subset\ZZ(\widehat\nu_{u_i}).
\]
Thus $\Lambda_i$ is orthogonal for $\nu_{u_i}$.
For each $i$, put
\[
 \eta_i=\mathop{*}_{u\ne u_i}\nu_u,\qquad \mu=\nu_{u_i}*\eta_i.
\]
All these measures are supported on $[0,\infty)$ and have zero in
their supports. Moreover, $\eta_i(\{0\})=0$: because $r>1$, its
random digit series has infinitely many nonnegative independent terms.
The sum can vanish only if every digit is zero. The probability that
the first $M$ such digits are zero is $N^{-M}$, which tends to zero.
These observations provide the endpoint conditions used below.

Take $f_\varepsilon=\mathbf1_{[0,\varepsilon]}$. 
Then $\|f_\varepsilon\|_2^2=\mu([0,\varepsilon])$ and $f_\varepsilon(x+y)>0$ when
$0\leq y\leq\varepsilon$ and $0\leq x\leq\varepsilon-y$.
For $\lambda\in\Lambda_i$, its Fourier coefficient in $L^2(\mu)$ is
\[
c_\lambda=\langle f_\varepsilon,e_\lambda\rangle_{L^2(\mu)}=
\int_{[0,\varepsilon]} e^{-2\pi\mathrm{i}\lambda y}
 \left(\int_{[0,\varepsilon-y]}e^{-2\pi\mathrm{i}\lambda x}
 \,d\nu_{u_i}(x)\right)d\eta_i(y).
\]
For each fixed $y\in[0,\varepsilon]$, the inner integrals are the
Fourier coefficients of $g_y=\mathbf1_{[0,\varepsilon-y]}$ in
$L^2(\nu_{u_i})$. Since $\Lambda_i$ is orthogonal for $\nu_{u_i}$, Bessel gives
\[
 \sum_{\lambda\in\Lambda_i}
 \left|\int_{[0,\varepsilon-y]}e^{-2\pi\mathrm{i}\lambda x}\,d\nu_{u_i}(x)\right|^2
 \leq\|g_y\|_{L^2(\nu_{u_i})}^2=\nu_{u_i}([0,\varepsilon-y]).
\]
Fix a finite subset $F\subset\Lambda_i$ and, for $0\leq y\leq\varepsilon$,
define
\[
 a_\lambda(y)=e^{-2\pi\mathrm{i}\lambda y}
 \int_{[0,\varepsilon-y]}e^{-2\pi\mathrm{i}\lambda x}\,d\nu_{u_i}(x),
 \qquad A_{i,F}(y)=(a_\lambda(y))_{\lambda\in F}.
\]
These functions are measurable and bounded, so their finite-dimensional
vector integral is defined coordinatewise. Since
$c_\lambda=\int_{[0,\varepsilon]}a_\lambda(y)\,d\eta_i(y)$,
Minkowski's integral inequality in $\ell^2(F)$ gives
\begin{align*}
 \left(\sum_{\lambda\in F}|c_\lambda|^2\right)^{1/2}
 &=\left\|\int_{[0,\varepsilon]}A_{i,F}(y)\,d\eta_i(y)
       \right\|_{\ell^2(F)}\\
 &\leq\int_{[0,\varepsilon]}\|A_{i,F}(y)\|_{\ell^2(F)}\,d\eta_i(y)\\
 &\leq\int_{[0,\varepsilon]}
       \nu_{u_i}([0,\varepsilon-y])^{1/2}\,d\eta_i(y).
\end{align*}
The last inequality follows from Bessel's inequality above, since the
phase factor in $a_\lambda(y)$ has modulus one. The bound is independent
of $F$. Taking the supremum over all finite subsets $F\subset\Lambda_i$
therefore yields
\[
 \left(\sum_{\lambda\in\Lambda_i}|c_\lambda|^2\right)^{1/2}
 \leq\int_{[0,\varepsilon]}
       \nu_{u_i}([0,\varepsilon-y])^{1/2}\,d\eta_i(y).
\]
This passage is valid because a sum of nonnegative terms is the
supremum of its finite subsums; no prior summability assumption on
$(c_\lambda)_{\lambda\in\Lambda_i}$ is needed.

By Cauchy--Schwarz and the convolution identity,
the square of the integral on the right is bounded by
\[
 \left(\int_{[0,\varepsilon]}1\,d\eta_i\right)
 \left(\int_{[0,\varepsilon]}\nu_{u_i}([0,\varepsilon-y])\,d\eta_i(y)\right)
 =\eta_i([0,\varepsilon])\mu([0,\varepsilon]).
\]
In particular,
\begin{equation}\label{(4.3)}
 \sum_{\lambda\in\Lambda_i}|c_\lambda|^2
 \leq\eta_i([0,\varepsilon])\mu([0,\varepsilon]).
\end{equation}
Since $s$ is finite and $\eta_i(\{0\})=0$, choose $\varepsilon>0$
such that $\sum_i\eta_i([0,\varepsilon])<1$.
Also $\mu([0,\varepsilon])>0$.
Here $\eta_i([0,\varepsilon])\to\eta_i(\{0\})=0$ by continuity
from above for finite measures. There are only finitely many $i$, so
their sum also tends to zero. Positivity of $\mu([0,\varepsilon])$
follows from $0\in\supp\mu$ and $\supp\mu\subset[0,\infty)$.
These two facts explain the strict inequality. Summing \eqref{(4.3)} gives
\[
 \sum_{\lambda\in\Lambda}|c_\lambda|^2
 <\mu([0,\varepsilon])=\|f_\varepsilon\|_{L^2(\mu)}^2,
\]
contradicting Parseval's identity if $\Lambda$ were a spectrum.
This contradiction excludes $r>1$, and therefore proves $r=1$.

\subsection{Excluding nonintegral rational reciprocals}
\begin{lem}\label{thm4.1}
Let $q>p\geq1$ be coprime integers, and suppose $0\in\Gamma$ with
\[
 (\Gamma-\Gamma)\setminus\{0\}
 \subset\bigcup_{n\geq1}(q/p)^n\mathbb Z.
\]
Then there is an integer $M\geq1$ such that
$\Gamma\subset p^{-M}\mathbb Z$.
\end{lem}
\begin{proof}
If $\Gamma=\{0\}$, take $M=1$. If $p=1$, the hypothesis gives
$\Gamma\subset\mathbb Z$, so again $M=1$ suffices. We may therefore
assume $p>1$ and $\Gamma\ne\{0\}$.
Each nonzero $\lambda\in\Gamma$ can be written
$\lambda=(q/p)^{n_\lambda}a_\lambda$, where $n_\lambda\geq1$,
$a_\lambda\in\mathbb Z$ and $q\nmid a_\lambda$.
To obtain this form, absorb all factors of $q$ in an initial coefficient
into the exponent, multiplying the new coefficient by the corresponding
power of $p$.
Let $M$ be the least $n_\lambda$, and fix
$\lambda_0=(q/p)^M a_0$.
For $n_\lambda>M$, write
\[
 (q/p)^{n_\lambda}a_\lambda-(q/p)^M a_0=(q/p)^k c,
 \qquad q\nmid c.
\]
We claim that $k=M$. Multiply the equation by $p^T$, where
$T=\max\{n_\lambda,k\}$, to obtain an equation of integers:
\[
 q^{n_\lambda}p^{T-n_\lambda}a_\lambda-q^Mp^{T-M}a_0
 =q^kp^{T-k}c.
\]
After division by $q^M$, the left side is not divisible by $q$:
its first term is divisible by $q$ because $n_\lambda>M$, whereas its
second is not because $q\nmid a_0$ and $\gcd(p,q)=1$.
Thus the largest power of $q$ dividing the left side of the integer
equation is exactly $q^M$. The corresponding power on the right is
exactly $q^k$, since $q\nmid c$. Hence $k=M$. 

Dividing the original difference equation by $(q/p)^M$ now gives
\[
 (q/p)^{n_\lambda-M}a_\lambda=a_0+c\in\mathbb Z.
\]
Consequently $\lambda=(q/p)^M(a_0+c)\in p^{-M}\mathbb Z$.
Elements represented with $n_\lambda=M$ plainly have the same
property, as does zero. This proves the lemma.
\end{proof}

\medskip\noindent\textbf{Finite partitions for a rational contraction.}
Let $\rho=p/q$ with $\gcd(p,q)=1$ and $p>1$.
Here $K=\mathbb Q$. Lemma~\ref{thm3.1} gives a finite partition
\begin{equation}\label{(4.4)}
 \Lambda=\bigcup_{i=1}^{s}\Lambda_i,
 \qquad\Lambda_i\subset\gamma_i+\mathbb Q,
\end{equation}
Choose the representatives so that $\gamma_1=0$. There is an integer
$H\geq1$ such that
\begin{equation}\label{(4.5)}
 \Lambda_i-\Lambda_j\subset\ZZ(\widehat\mu_H)
 \qquad(i\ne j).
\end{equation}
If $s=1$, choose $H=1$.
If $Z\cap\mathbb Q$ is empty, every class in \eqref{(4.4)} has at most
one element, contradicting infinitude of $\Lambda$. Otherwise choose a
common denominator $L\geq2$ of $Z\cap\mathbb Q$.

For $\Gamma_i=\Lambda_i-\gamma_i\subset\mathbb Q$, every nonzero
difference of two elements of $\Gamma_i$ has the form
\begin{equation}\label{(4.6)}
 (q/p)^n(z+m),\qquad n\geq1,\quad z\in Z\cap\mathbb Q,\quad m\in\mathbb Z.
\end{equation}
Apply Lemma~\ref{thm4.1} to each $L\Gamma_i$. Taking the maximum of
the finitely many denominator exponents gives a common integer
$M\geq1$ such that
\begin{equation}\label{(4.7)}
 \Gamma_i\subset (Lp^M)^{-1}\mathbb Z
\end{equation}
for any $1\leq i\leq s$. 
For integers $n\geq1$, define the finite sets
\begin{equation}\label{(4.8)}
 B_n=\bigcup_{i=1}^{s}
 \{\lambda\in\Lambda_i:|\lambda-\gamma_i|\leq q^n\}.
\end{equation}
They increase to $\Lambda$ and contain zero, as $\gamma_1=0$.

\medskip\noindent\textbf{The finite-section bound.}
We claim that there is a fixed integer $C_0\geq H$ such that $B_n$ is
orthogonal for $\mu_{n+C_0}$ for every $n\geq1$. 
Take distinct $\lambda,\lambda'\in B_n$, with $\lambda\in\Lambda_i$
and $\lambda'\in\Lambda_{i'}$. If $i\ne i'$, \eqref{(4.5)} gives
\begin{align*}
\lambda'-\lambda\in\mathcal{Z}(\widehat{\mu_H})\subset\mathcal{Z}(\widehat{\mu_{n+C_0}}).
\end{align*}
Suppose next that $i=i'$. 
Put $\gamma=\lambda-\lambda'$. Then
\[
 0<|\gamma|\leq2q^n,\qquad Lp^M\gamma\in\mathbb Z.
\]
From \eqref{(4.6)}, $\gamma$ can be written
$(q/p)^j a/L$ with $a\in\mathbb Z\setminus\{0\}$.
Thus
\begin{equation}\label{(4.9)}
 Lp^M\gamma=q^j p^{M-j}a\in q^j\mathbb Z.
\end{equation}
The divisibility in \eqref{(4.9)} requires a check when $j>M$.
Put $t=Lp^M\gamma\in\mathbb Z\setminus\{0\}$. If $j\leq M$,
the divisibility is immediate. If $j>M$, then
$p^{j-M}t=q^ja$, and coprimality of $p$ and $q$ implies $q^j\mid t$.
Thus $t$ is a nonzero integer multiple of $q^j$, and consequently
\[
 q^j\leq|Lp^M\gamma|\leq 2Lp^M q^n,
 \qquad j\leq n+\log_q(2Lp^M).
\]
Taking
\begin{equation}\label{(4.10)}
 C_0\geq\max\{H,1,\lceil\log_q(2Lp^M)\rceil\}
\end{equation}
proves the claim. 

\medskip\noindent\textbf{Fourier decay and failure of completeness.}
We use the following two analytic inputs, in precisely the form needed
for the rational-contraction argument.
\begin{lem}[{\cite{GM17}}]\label{thm4.2}
Let $\rho=p/q\in(0,1)$ be in lowest terms with $p>1$, and let
$D\subset\mathbb Z$ contain at least two elements. There are constants
$C,\gamma_0>0$ such that
\begin{equation}\label{(4.11)}
 |\widehat\mu_{\rho,D}(\xi)|
 \leq C\bigl(\log(2+|\xi|)\bigr)^{-\gamma_0}\qquad(\xi\in\mathbb R).
\end{equation}
\end{lem}
We quote this logarithmic Fourier-decay estimate directly from
\cite{GM17}. The negative exponent in \eqref{(4.11)} means that the
bound tends to zero as $|\xi|\to\infty$.

\begin{lem}[{\cite[Theorem~2.5(ii)]{AW21}}]\label{thm4.3}
Let $\Lambda$ be an orthogonal set for $\mu=\mu_{\rho,D}$.
Suppose $n_k$ are increasing positive integers and $A_k$ are increasing
finite sets with $0\in A_k$, $\bigcup_k A_k=\Lambda$, and $A_k$
orthogonal for $\mu_{n_k}$. If, for all sufficiently large $k$,
\begin{equation}\label{(4.12)}
 \sup_{\lambda\in A_{k+1}\setminus A_k}
 |\widehat\mu(\rho^{n_{k+1}}\lambda)|^2\leq c_k,
 \qquad \sum_k c_k<\infty,
\end{equation}
then $\Lambda$ is not a spectrum of $\mu$.
The supremum over an empty difference is interpreted as zero.
\end{lem}
This follows from \cite[Theorem~2.5(ii)]{AW21}, attributed there to
\cite{AHH19}. To match its notation, set $\widetilde\Lambda_n=\{0\}$
for $n<n_1$ and $\widetilde\Lambda_n=A_k$ when
$n_k\leq n<n_{k+1}$. Since
$\ZZ(\widehat\mu_{n_k})\subset\ZZ(\widehat\mu_n)$ for $n\geq n_k$,
each $\widetilde\Lambda_n$ is orthogonal for $\mu_n$.
These sets are nested and exhaust $\Lambda$. Applying the cited
criterion along the subsequence $n_k$ gives the stated form.
We now verify its assumptions for our sets.

Let $\gamma_0>0$ be as in Lemma~\ref{thm4.2}. Choose an integer $a\geq2$
with $2a\gamma_0>1$, and set
\begin{equation}\label{(4.13)}
 n_k=k^a+C_0,\qquad A_k=B_{k^a}.
\end{equation}
The previous subsection shows that $A_k$ is orthogonal for $\mu_{n_k}$.
Let $G=\max_{1\leq i\leq s}|\gamma_i|$. If $\lambda\in A_{k+1}\setminus A_k$,
its class in \eqref{(4.4)} is unique, so
\[
 |\lambda|\geq q^{k^a}-G\geq\tfrac12 q^{k^a}
\]
for all sufficiently large $k$. Consequently
\begin{align}
 \log\bigl(\rho^{n_{k+1}}|\lambda|\bigr)
 &\geq ((k+1)^a+C_0)\log(p/q)+k^a\log q-\log2\notag\\
 &= (k+1)^a\log p-\bigl((k+1)^a-k^a\bigr)\log q+O(1)\notag\\
 &\geq c k^a                                                     \label{(4.14)}
\end{align}
with $c>0$, because $\log p>0$ and
$(k+1)^a-k^a=O(k^{a-1})$.
By Lemma~\ref{thm4.2},
\begin{equation}\label{(4.15)}
 \sup_{\lambda\in A_{k+1}\setminus A_k}
 |\widehat\mu(\rho^{n_{k+1}}\lambda)|^2
 \leq C_1 k^{-2a\gamma_0}.
\end{equation}
The right-hand side is summable because $2a\gamma_0>1$.

We have now checked all assumptions of Lemma~\ref{thm4.3}: the levels
$n_k$ increase, the finite sets $A_k$ increase and exhaust $\Lambda$,
they contain zero, they are orthogonal at the required finite levels,
and the tail bounds are summable. Any finitely many initial values of
$k$ are irrelevant to that criterion. Lemma~\ref{thm4.3} implies that
$\Lambda$ is not a spectrum, a contradiction.
Hence $p=1$. Together with $r=1$ in \eqref{(4.1)}, this gives
$\rho^{-1}=q\in\mathbb N$, $q\geq2$, and completes the proof of
Theorem~\ref{thm1.1}.\hfill$\square$

\section*{Declarations}
\subsection*{AI use}

The author used ChatGPT (OpenAI) to assist with manuscript
organization, language editing, and the completion of one
part of the proof. Before using this tool, the author had
established the spectral decomposition lemma(Lemma \ref{thm3.1}),
Theorem~\ref{thm1.2}, and the case $\rho=p/q$ with $p>1$.
For the remaining case, where $\rho\notin\mathbb{Q}$ but
$\rho^r\in\mathbb{Q}$ for some integer $r\geq2$, ChatGPT
supplied the missing proof details based on the author's
spectral decomposition outline. The author checked this
proof in detail and verified its correctness, and takes
full responsibility for the final manuscript.
\subsection*{Conflict of Interest}
The author declares no conflict of interest.

\end{document}